\documentclass[psamsfonts]{amsart}

\usepackage{amssymb,amsfonts}
\usepackage[all,arc]{xy}
\usepackage{enumerate}
\usepackage{mathrsfs}
\usepackage{amsmath}

\theoremstyle{definition}

\theoremstyle{remark}

\makeatletter
\let\c@equation\c@thm
\makeatother
\numberwithin{equation}{section}

\title{A Polynomial Method Approach to Zero-Sum Subsets in ${\mathbb{F}}_{p}^{2}$}

\author{Cosmin Pohoata}

\begin{document}

\begin{abstract}
In this paper we prove that every subset of
$\mathbb{F}_p^2$ meeting all $p+1$ lines passing through the origin has a zero-sum subset. This is motivated by a result of Gao, Ruzsa and Thangadurai which states that $OL(\mathbb{F}_{p}^{2})=p+OL(\mathbb{F}_{p})-1$, for sufficiently large primes $p$. Here $OL(G)$ denotes the so-called Olson constant of the additive group $G$ and represents the smallest integer such that no subset of cardinality $OL(G)$ is zero-sum-free. Our proof is in the spirit of the Combinatorial Nullstellensatz.
\end{abstract}

\maketitle

\section{Introduction}

\bigskip

\bigskip

For a set $A$ in an additive group $G$, denote 
$$\Sigma(A)=\left\{\sum_{x\in B}x\ |\ B\subset A\right\},\ \Sigma^*(A)=\left\{\sum_{x\in B}x\ |\ B\subset A, B\ne \emptyset\right\}.$$
Define the \emph{Olson constant} $OL(G)$ of a finite additive group $G$ to be the minimal $d$ such that every subset $A\subset G$ of cardinality $d$ satisfies $0\in \Sigma^*(A)$. 

It turns out this is a difficult quantity to compute for general groups. First, Erd\H{o}s and Heilbronn initially proved in [5] that there exists an absolute constant $c$ such that $OL(\mathbb{F}_{p}) \leq c\sqrt{p}$, where $p$ is a prime. They also conjectured that this bound extends to arbitrary additive groups $G$, and years later Szemer\' edi showed in [12] that $OL(G) \leq k\sqrt{|G|}$ for some $k>0$. However, both results were far from being optimal. Olson [11] also came into the picture later and proved that that $OL(G) \leq 2\sqrt{|G|}$, which has been subsequently improved by Hamidoune and Zemor [7], who proved that
$$OL(G) \leq \sqrt{2|G|} + O(|G|^{1/3} \log|G|),$$
conjecturing that $c=\sqrt{2}$ is perhaps optimal.

The story is not much better for groups of the form $\mathbb{F}_{p}^{d}$, where $p$ is an odd prime. For $d=1$, Hamidoune and Zemor were indeed correct; in [8], Nguyen, Szemeredi and Vu improved upon the Erd\H{o}s-Heilbronn theorem and showed that for sufficiently large primes $p$, $OL(\mathbb{F}_{p}) \leq \sqrt{2p}$. This was subsequently also shown by Balandraud in [2] for all primes $p$, who in fact actually proved that for $A\subset \mathbb{F}_p$, with $|A|=d$ and such that $A\cap (-A)=\emptyset$, we have that
$$|\Sigma(A)|\geq \min\left\{p,d(d+1)/2+1\right\}\ \ \text{and}\ \ |\Sigma^*(A)|\geq \min\left\{p,d(d+1)/2\right\}.$$
Balandraud's original proof made use of some impressive recursions for binomial determinants, but it has been recently recasted by the same author in [3] as a natural consequence of the Combinatorial Nullstellensatz due to Alon [2]. For $d=2$, Gao, Ruzsa and Thangadurai showed in [7] that 
$$OL(\mathbb{F}_{p}^{2})=p+OL(\mathbb{F}_{p})-1$$
holds for $p > 4.67 \times 10^{34}$. Their argument was subsequently updated by Bhowmik and Schlage-Puchta to show that the equality holds for all primes $>6000$; thus, given the result for $\mathbb{F}_{p}$, the Olson constant constant of $\mathbb{F}_{p}^{2}$ is also determined for large primes. For higher dimensions, Gao, Ruzsa and Thangadurai conjectured that 
$$OL(\mathbb{F}_{p}^{d}) = p + OL(\mathbb{F}_{p}^{d-1})-1$$
should hold in general, but little is known in this case. We refer to [9] for an account of the state of the art. 

In this paper we shall prove the following

\bigskip

{\bf{Theorem 1}}. {\it{Let $p$ be an odd prime and let $\left\{x_{1},\ldots,x_{n}\right\} \subset \mathbb{F}_{p}^{2}\ \backslash \left\{0,0\right\}$ be a set of vectors spanning all $p+1$ nonzero directions from $0$ in $\mathbb{F}_{p}^{2}$. We prove that there exists a subset $I \subset \left\{1,\ldots, n\right\}$ such that
$$\sum_{i \in I}{x_{i}}=0.$$}}

\bigskip

This is meant to complement the Gao-Rusza-Thangadura-Bhowmik-Schlage-Puchta theorem. Here, we shall say that two vectors $(a,b), (c,d) \in \mathbb{F}_{p}^{2}$ have the same direction from zero if there is some $t \in \mathbb{F}_{p}\ \backslash \left\{0\right\}$ such that $(a,b) = (tc,td)$. The above result thus says that a set $A\subset \mathbb{F}_p^2$, with $(0,0)\notin A$ and with the property that it meets all $p+1$ lines of the form $ax+by=0$, must contain a non-empty subset $B\subset A$ with zero sum.

\bigskip

\section{Preliminaries}

\bigskip

Although not implicitly used in our proof, we first recall the statement of the Combinatorial Nullstellensatz.

\bigskip

{\bf{Theorem 2}}. {\it{Let $\mathbb{F}$ be an arbitrary field, and let $f = f(x_{1},\ldots,x_{n})$ be a polynomial in $F[x_{1},\ldots,x_{n}]$. Suppose the degree $\deg(f)$ of $f$ is $\sum_{i=1}^{n}{t_{i}}$, where each $t_{i}$ is a nonnegative integer, and suppose the coefficient of $\prod_{i=1}^{n}{x_{i\
}^{t_{i}}}$ in $f$ is nonzero. Then, if $S_{1},\ldots,S_{n}$ are the subsets of $F$ with $|S_{i}| > t_{i}$, there are $s_{1} \ldots S_{1}$,$\ldots$, $s_{n} \in S_{n}$ so that
$$f(s_{1},\ldots,s_{n}) \neq 0.$$}}

\bigskip

This is usually proven in general by induction on the number of variables (see for instance [1]), but whenever $S_{i} = \left\{0,1\right\}$, $t_{i}=1$ for each $i \in \left\{1,\ldots,n\right\}$, there's a very neat way of thinking about this, which will be relevant to our proof of Theorem 1. One can simply consider the following ``integral'' of $f$ on $\left\{0,1\right\}^{n}$:
$$\mathbb{I}_{f} := \sum_{(x_1,\ldots,x_{n})\in \left\{0,1\right\}^{n}}{(-1)^{x_1+\dots+x_{n}} f(x_{1},\ldots,x_{n})}.$$
If $f$ vanishes everywhere on $S_{1} \times \ldots \times S_{n}$, then clearly $\mathbb{I}_{f} = 0$. 

On the other hand, if $c$ is the nonzero coefficient of $\prod_{i=1}^{n}{x_{i}}$ in $f$, then
$$g(x_{1},\ldots,x_{n}):= f(x_{1},\ldots,x_{n}) - c\prod_{i=1}^{n} x_{i}$$
is a polynomial of degree less than $n$ in $\mathbb{F}[x_{1},\ldots,x_{n}]$. In particular, every monomial appearing in $g \in \mathbb{F}[x_{1},\ldots,x_{n}]$ with nonzero coefficient has to be of the form $x_{1}^{i_{1}}\ldots x_{n}^{i_{n}}$, with $i_{j} \in \left\{0,1\right\}$ not all $1$. For such a monomial, where say $i_{K}=0$ for some $K \in \left\{1,\ldots,n\right\}$, we then note that
\begin{eqnarray*}
&& \sum_{(x_1,\ldots,x_{n})\in \{0,1\}^{n}} {(-1)^{x_1+\dots+x_{n}} x_{1}^{i_{1}}\ldots x_{n}^{i_{n}}}\\
&=& \sum_{x_{K} \in \left\{0,1\right\} } (-1)^{x_{K}} \sum_{x_{j} \in \left\{0,1\right\}, j \neq K} (-1)^{\sum_{j \neq K}x_{j} }x_{1}^{i_{1}}\ldots x_{n}^{i_{n}}\\
&=& 0.
\end{eqnarray*}
Therefore, $\mathbb{I}_{g} = 0$, which yields $\mathbb{I}_{f} = \mathbb{I}_{c\prod_{i=1}^{n}{x_{i}}} = \pm c \neq 0$, a contradiction.

\bigskip

We will consider such an integral in our proof of Theorem 1 as well. Before doing so however, let us mention first that one of the first certificates of the usefulness of the Combinatorial Nullstellensatz in [1] was a proof of the following theorem due to Chevalley and Warning.

\bigskip

{\bf{Theorem 3}}. {\it{Let $p$ be a prime and let $\mathbb{F}$ be a finite field. Let $$P_{1}=P_{1}(x_{1},\ldots,x_{n}),\ldots,P_{m}=P_{m}(x_{1},\ldots,x_{n})$$
be $m$ polynomials in $\mathbb{F}[x_{1},\ldots,x_{n}]$ such that
$$n > \sum_{j=1}^{m}{d_{j}},$$
where $d_{j}$ is the total degree of $f_{j}$. Then, if the polynomials $P_{i}$ have a common zero $(c_{1},\ldots,c_{n})$, they must have another common zero.}}

\bigskip

We mention this because the one can use the Theorem 3 as a ``black-box'' to show something in the spirit of Theorem 1.

\bigskip

{\bf{Theorem 4}}. {\it{Let $p$ be a prime and let $n$ be integers such that $n > 2(p-1)$. Then, among any $n$ elements $v_{1},\ldots,v_{n}$ of $\mathbb{F}_{p}^{2}$ there exists a a nonempty subsequence with a zero-sum.}}

\bigskip

In fact, it is a result of Olson [10] that if $G$ is an abelian $p$-group of the form $G = \mathbb{Z}_{p^{\alpha_{1}}} \oplus \ldots \oplus \mathbb{Z}_{p^{\alpha_{k}}}$, and $n$ is any integer such that 
$$n > 1 + \sum_{i=1}^{k}{(p^{\alpha_{i}}-1)},$$
then among every sequence of elements of $G$ with length $n$ there is a nonempty subsequence with a zero sum. We sketch the proof for the scenario described in Theorem 4. 

\bigskip

{\it{Proof of Theorem 4}}. Assume without loss of generality that $n = 2p-1$, and furthermore let $v_{i} = (a_{i1},a_{i2})$, for each $1 \leq i \leq n$. Define
$$f_{j}(x_{1},\ldots,x_{n}):=\sum_{i=1}^{n}{a_{ij}x_{i}^{p-1}}$$
for each $j \in \left\{1,2\right\}$. Note that these polynomials fulfll the hypothesis of the Chevalley-Warning theorem: the sum of their degrees is $2(p-1)$, which is less than $n = 2p-1$, the number of variables. Also, note that the system $f_{1}=f_{2}=0$ has the trivial solution $x_{1}=\ldots=x_{n}=0$. By the Chevalley-Warning theorem, it thus has another one, say $(\alpha_{1},\ldots,\alpha_{n})$. Setting $I:=\left\{i\ |\ \alpha_{i} \neq 0\right\},$
this is non-empty and it is easy to check that
$$\sum_{i \in I}{v_{i}} = 0,$$
as desired.

$\hfill \square$

\bigskip

\bigskip

\section{Proof of Theorem 1}

\bigskip

Let $A$ be a set of vectors $\left\{x_{1},\ldots,x_{n}\right\} \subset \mathbb{F}_{p}^{2}\ \backslash \left\{0,0\right\}$ spanning all $p+1$ nonzero directions from $0$ in $\mathbb{F}_{p}^{2}$. Recall that we want to prove that there exists a subset $I \subset \left\{1,\ldots, n\right\}$ such that
$$\sum_{i \in I}{x_{i}}=0.$$

Without loss of generality, assume that $n=p+1$ and that 
$$A = \left\{(a_{1},a_{1}),\ldots,(a_{p-1},(p-1)a_{p-1}),(a_{p},0),(0,a_{p+1})\right\}$$
where $a_{1},\ldots,a_{p+1}$ are non-zero elements in $\mathbb{F}_{p}$. We claim that there exists a non-empty subset of $A$ with sum of elements equal to zero.

\smallskip

We will first prove a crucial Lemma.

\bigskip

{\bf{Lemma 5}}. {\it{The polynomial $f \in \mathbb{F}_{p}[x_{1},\ldots,x_{p+1}]$ defined by
$$f(x_{1},\ldots,x_{p+1}) = (x_{1}+\ldots+x_{p})^{p-1} (x_{1}+2x_{2}+\ldots+px_{p}+x_{p+1})^{p-1}$$
contains no monomial where variables $x_{1},\ldots,x_{p+1}$ all have positive exponents.}}

\bigskip

{\it{Proof of Lemma 5}}. Assume that the coefficient of some monomial $x_1^{c_1}\dots x_{p+1}^{c_{p+1}}$ appearing in the expansion of $f$ is nonzero modulo $p$, where $c_i>0$ for each $i \in \left\{1,\ldots,p+1\right\}$. First, note that each $c_i$ is also less than $p$, since otherwise the monomial would have degree at least $2p > 2p-2 = \deg f$. Moreover, the coefficient $[x_1^{c_1}\dots x_{p+1}^{c_{p+1}}](f)$ of $x_1^{c_1}\dots x_{p+1}^{c_{p+1}}$ in $f$ satisfies
$$[x_1^{c_1}\dots x_{p+1}^{c_{p+1}}](f) = \left(c_{1}\cdot \ldots \cdot c_{p+1}\right)^{-1} [x_1^{c_1-1}\dots x_{p+1}^{c_{p+1}-1}](g),$$ where the latter denotes the coefficient of $x_1^{c_1-1}\dots x_{p+1}^{c_{p+1}-1}$ in the derivative
$$g(x_{1},\ldots,x_{p+1}):= \frac{\partial^{n} f(x_{1},\ldots,x_{p+1})}{\partial x_{1}\ldots \partial x_{p+1}}.$$
Since $p > c_{i} > 0$ for each $i$, this coefficient is also nonzero modulo $p$. 

On the other hand, the polynomial $g$ is identically zero. Indeed, in the computation of the multiple derivative, we differentiate $(x_1+x_2+\dots+x_{p})^{p-1}$ with respect to some $k$ variables ($2\leq k\leq p-1$), including $x_p$ but excluding $x_{p+1}$, and we differentiate $(x_1+2x_2+\dots+px_p+x_{p+1})^{p-1}$ with respect to other $p+1-k$ variables (including $x_{p+1}$ but excluding $x_{p}$). It is easy to check that for each summand appearing in $g$, a multiple of the $(p-k)$-th elementary symmetric polynomials in the variables $1,2,\dots,p$ arises for some $k \in \left\{2,\ldots,p-1\right\}$. By Vieta's formulas, these are all equal to $0$ over $\mathbb{F}_{p}$, as $(t-1)\dots(t-p)=t^p-t$. We have therefore arrived at a contradiction.

$\hfill \square$

\bigskip

Now, suppose that $A$ contains no non-empty subset with sum of elements equal to zero. Consider the following polynomial in $\mathbb{F}_{p}[x_{1},\ldots,x_{p+1}]$:
$$
P(x_1,\dots,x_{p+1})=\left[1- \left(\sum_{i=1}^{p}{a_{i}x_{i}}\right)^{p-1}\right]\left[1-\left(\sum_{j=1}^{p+1}{ja_{j}x_{j}}\right)^{p-1}\right].$$
Note that $P(x_{1},\ldots,x_{p+1})=1$ precisely when $\sum_{i=1}^{p}{a_{i}x_{i}}=\sum_{j=1}^{p+1}{ja_{j}x_{j}}=0$ and is $0$ otherwise; therefore, when restricted on $\left\{0,1\right\}^{p+1}$, the polynomial $P$ identifies zero-sum subsets of $A$. In particular, given our zero-sum-subset free assumption on $A$, we have that 
$$P(0,\ldots,0)=1\ \ \ \text{and}\ \ \ P(x_{1},\ldots,x_{p+1})=0,$$
for every $(p+1)$-tuple $(x_{1},\ldots,x_{p+1}) \in \left\{0,1\right\}^{p+1} \backslash \left\{0,\ldots,0\right\}$. 
In particular,
$$\sum_{(x_1,\ldots,x_{p+1})\in \{0,1\}^{p+1}}{(-1)^{x_1+\dots+x_{p+1}} P(x_1,\dots,x_{p+1})}=1.$$

If $Q \in \mathbb{F}_{p}[x_{1},\ldots,x_{p+1}]$ is defined by
$$
Q(x_1,\dots,x_{p+1})=\left(\sum_{i=1}^{p}{a_i x_{i}}\right)^{p-1}\left(\sum_{j=1}^{p+1}{ja_{j} x_{j}}\right)^{p-1},$$
then $P-Q$ is a polynomial of degree $p-1$ in $\mathbb{F}_{p}[x_{1},\ldots,x_{p+1}]$. In particular, there is no monomial in $P-Q$ which contains all variables $x_{1},\ldots,x_{p+1}$, each with positive exponent. For each monomial $x_{1}^{i_{1}}\ldots x_{p+1}^{i_{p+1}}$, for which $i_{K}=0$ for some $K \in \left\{1,\ldots,p+1\right\}$, we then note that
\begin{eqnarray*}
&& \sum_{(x_1,\ldots,x_{p+1})\in \{0,1\}^{p+1}}{(-1)^{x_1+\dots+x_{p+1}} x_{1}^{i_{1}}\ldots x_{p+1}^{i_{p+1}}}\\
&=& \sum_{x_{K} \in \left\{0,1\right\}} (-1)^{x_{K}} \sum_{x_{j} \in \left\{0,1\right\}, j \neq K}{(-1)^{x_1+\dots+x_{p+1}} x_{1}^{i_{1}}\ldots x_{p+1}^{i_{p+1}}}\\
&=& 0.
\end{eqnarray*}
Since this holds for every monomial appearing in $P-Q$, we get that
$$\sum_{(x_1,\ldots,x_{p+1})\in \{0,1\}^{p+1}}{(-1)^{x_1+\dots+x_{p+1}} \left(P(x_1,\dots,x_{p+1}) - Q(x_1,\dots,x_{p+1})\right)}=0.$$
Consequently,
$$\sum_{(x_1,\ldots,x_{p+1})\in \{0,1\}^{p+1}}{(-1)^{x_1+\dots+x_{p+1}} Q(x_1,\dots,x_{p+1})}=1.$$
But this gives a contradiction with the structure of $A$. Lemma 5 says that each monomial appearing in $Q$ must have at least one variable absent; by performing the same trick as above, we can thus argue that for each monomial $x_{1}^{j_{1}}\ldots x_{p+1}^{j_{p+1}}$ appearing in $Q$, we have that
$$\sum_{(x_1,\ldots,x_{p+1})\in \{0,1\}^{p+1}}{(-1)^{x_1+\dots+x_{p+1}} x_{1}^{j_{1}}\ldots x_{p+1}^{j_{p+1}}}=0,$$
and therefore
$$\sum_{(x_1,\ldots,x_{p+1})\in \{0,1\}^{p+1}}{(-1)^{x_1+\dots+x_{p+1}} Q(x_1,\dots,x_{p+1})}=0,$$
which is a contradiction. 

$\hfill \square$

\bigskip

\bigskip

\subsection*{Acknowledgments} 

\bigskip

I would like to thank Fedor Petrov for helpful comments on a prior version of this preprint.

\bigskip

\bigskip

\small{ \textsc{California Institute of Technology, Pasadena, CA}

{\it{E-mail address}}: apohoata@caltech.edu}

\end{document}